\documentclass{article}

\usepackage[top=2.5cm, bottom=2cm, inner=2.5cm, outer=2.5cm]{geometry}
\usepackage{color}
\usepackage{paralist,parskip}
\usepackage{booktabs, multirow, makecell}
\usepackage{subfigure, bm, comment}
\usepackage{amsmath,amsfonts,amssymb}
\usepackage{mathabx,fixmath,algorithm, algorithmic}
\usepackage{pstool}
\usepackage{hyperref}
\usepackage{cite}
\def\x{\mathbold{x}}
\def\A{\mathbold{A}}

\def\Q{\mathbf{Q}}

\def\g{\mathbold{g}}

\def\Z{\mathbf{Z}}

\newcommand{\w}{\mathbold{w}}

\usepackage{theorem}

\newtheorem{assumption}{Assumption}

\newtheorem{theorem}{Theorem}

\newtheorem{corollary}{Corollary}

\usepackage{enumitem}

\title{Time-Varying Optimization:\\ Algorithms and Engineering Applications}
\author{Andrea Simonetto, \emph{IBM Research Ireland, Dublin}}
\begin{document}

\setlength{\parskip}{2mm} 
\setlength{\parindent}{0pt}

\maketitle

\begin{abstract}
This is the write-up of the talk I gave at the 23rd International Symposium on Mathematical Programming (ISMP) in Bordeaux, France, July 6th, 2018. The talk was a general overview of the state of the art of time-varying, mainly convex, optimization, with special emphasis on discrete-time algorithms and applications in energy and transportation. This write-up is mathematically correct, while its style is somewhat less formal than a standard paper. 
\end{abstract}

\paragraph{Acknowledgements.} Many thanks to all the collaborators that are or have been working with me on these themes. Their names appear in the references at the end of this write-up. 

\section{Introduction}

With time-varying optimization, we mean the task of finding the minumum of an optimization problem that changes continuously in time. Let $f: \mathbf{R}^n \times \mathbf{R}_{+} \to \mathbf{R}$ be a convex function parametrized over time, i.e., $f(x; t)$, where $x \in \mathbf{R}^n$ is the decision variable and $t\geq 0$ is time. Let $X(t)\subseteq \mathbf{R}^n$ be a convex set, also changing in time. Then the problem at hand can be formulated as finding
\begin{equation}\label{tvp}
\min_{x\in X(t)} \, f(x; t), \quad \textrm{for all } t\geq 0.
\end{equation}
That is, we want to find the minimum at each point in time. These types of problems appear naturally in many applications, for example energy, robotics, transportation, as we will see. 

In this talk, we will sample Problem~\eqref{tvp} at defined sampling times $t_k$, with $k = 0, 1, \ldots$ and sampling period $h = t_{k+1} - t_k$, and arrive at a sequence of time-invariant problems
\begin{equation}\label{tip}
\min_{x\in X_k} \, f(x; t_k).
\end{equation}

When one can sample Problem~\eqref{tvp} at the desired sampling frequency and solve the resulting time-invariant problems~\eqref{tip} at the desired accuracy within the sampling period, we are in a batch solution mode. This batch approach is hardly viable, except for low dimensional problems that can be sampled with sufficient long sampling periods (i.e., when the problem changes sufficiently slowly). We won't follow this approach, instead we will pursue an on-line approach, which will find approximate solutions of each of the time-invariant problems~\eqref{tip} and eventually will get close to the minimum trajectory. 

To generate approximate solutions we will use running (or correction-only, or catching-up) algorithms and prediction-correction algorithms. We will touch upon primal and dual algorithms.

\subsection{Background}

Time-varying optimization has been around for quite some time, e.g.,~\cite{Polyak1987}. 

Continuous-time platforms have been discussed, e.g., in~\cite{Ye2015,Rahili2015a,Rahili2015,Gong2016,Fazlyab2015,Fazlyab2016}.

Running methods on discrete-time platforms can be traced back to Moreau~\cite{Moreau1977}, and subsequently have appeared in many contexts~\cite{Popkov2005,Tu2011,Bajovic2011, Dontchev2013,Zavlanos2013,Jakubiec2013,Ling2013,Simonetto2014c,Simonetto2014d,Ye2015,Xi2016a,Sun2017,Maros2017}.

A recent and fairly complete treatment is in~\cite{Simonetto20XX} (which forms also the basis for part of the results I will present here).

Prediction-correction methods arise from non-stationary optimization~\cite{Polyak1987,Popkov2005}, parametric programming~\cite{Robinson1980,Dontchev2009,Zavala2010,Dontchev2013,Kungurtsev2017}, and continuation methods in numerical mathematics~\cite{Allgower1990}. It also resembles evolutionary variational inequalities~\cite{Cojocaru2005,Nagurney2006} and path-following methods in interior point solvers~\cite{Nesterov2012}.

Part of the work that I present here is in~\cite{Paper1,Paper2,Paper3,Paper4}.

\section{Formulation}

Our main starting point is the time-invariant problem~\eqref{tip}. Unless otherwise said, we will assume that the function is nicely behaving, that for us means that 

\begin{assumption}\label{as.1}
Function $f(x; t)$ is $m$ strongly convex ($m>0$) and $L$ strongly smooth ($L>0$) over $x \in \mathbf{R}^n$, uniformly in time. 
\end{assumption}

Assumption~\ref{as.1} guarantees that the solution (i.e., the minimizer of~\eqref{tip}) exists and its unique at every time $t_k$ (of course assuming that $X_k$ is non-empty). This implies also that the \emph{solution trajectory} $x^*(t)$ of~\eqref{tvp} is well-defined. 

To see situations for which this is not true see~\cite{Guddat1990}.

\section{Running algorithms}

Running algorithms or, as we said, correction-only/catching-up algorithms have always the same prototype structure. Here I present the work of~\cite{Simonetto20XX} (so theorems and results are properly defined there); see the original paper for references to previous work.   

Running algorithms start with a approximate solution $x_0$ and generate a sequence $\{x_k\}$ by acquiring a new function at time $t_{k+1}$ and performing $C$ iterations of the selected method. For example, for the case of the running projected gradient, one does the following

\begin{itemize}
\item Time $t_0$, guess $x_0$
\item Time $t_{k+1}$
\begin{enumerate}
\item Acquire a new function $f(\cdot; t_{k+1})$ and the constraint set $X_{k+1}$
\item Set $y_{0} = x_{k}$
\item Perform $C$ times:
\begin{equation}\label{run}
y_{i+1} = \Pi_{X_{k+1}}[y_i - \alpha \nabla_{x} f(y_i; t_{k+1})]
\end{equation}
\item Set $x_{k+1} = y_{C}$
\end{enumerate}
\end{itemize}

In~\eqref{run}, $\alpha>0$ is the stepsize, while $\Pi_{X}$ is the projection onto the convex set $X$. 

\subsection{Theoretical results}

Typical theoretical results of running algorithms go as follow. Assume that the optimizer trajectory is well-behaved, that is that

\begin{assumption}\label{as.2}
The change in the optimizers of~\eqref{tip} is upper bounded as
$$
\|x^*(t_{k+1}) - x^*(t_k)\| \leq K, \quad \forall k\geq0.
$$ 
\end{assumption}

Assumption~\ref{as.2} guarantees that the optimizers are indeed trackable. Note that $K$ can be big, so we are not limited to small variations; on the other hand, the bounds will be big too if $K$ is big. 

Then, we have

\begin{theorem}\label{th.1}\emph{(Informal)} If your favorite method $\mathcal{M}$ converges Q-linearly to the optimizer of a time-invariant problem as
$$
\|x_k - x^*(t_k)\| \leq \varrho^C \|x_{k-1} - x^*(t_k)\|, \quad \varrho <1, 
$$
then the same method $\mathcal{M}$ converges Q-linearly to the optimizer trajectory of a time-varying problem up to an error bound as
$$
\|x_k - x^*(t_k)\| \leq \varrho^C (\|x_{k-1} - x^*(t_{k-1})\|+K), 
$$
and
$$
\limsup_{k\to\infty} \|x_k - x^*(t_k)\| = \varrho^C O(K).
$$
\end{theorem}

The theorem is fairly general and its based on the triangle inequality. What it says is that the sequence $\{x_k\}$ will track the solution trajectory up to a ball of size $\varrho^C O(K)$. If $C \to \infty$, we solve the time-invariant problem exactly and we are back to the time-invariant/batch mode (and the error is $0$). 

Based on this theorem, one can derive a corollary for the projected gradient method

\begin{corollary}
For $\alpha < 2/L$, the projected gradient method applied in a running mode generate a sequence $\{x_k\}$ that converges to the error bound $\varrho^C O(K)$ Q-linearly, with rate $\varrho = \max\{|1-\alpha m|,|1-\alpha L|\}$.  
\end{corollary}

Equation~\eqref{run} can be substituted with other methods, and Theorem~\ref{th.1} is true for a variety of methods $\mathcal{M}$, such as
\begin{itemize}
\item Proximal point method, for $f(\cdot; t_{k})$ strongly convex, $X_{k} \subseteq \mathbf{R}^n$;

\item Forward-backward splitting (minimizing $f(x; t) + g(x; t)$), for $f(\cdot; t_{k})$ strongly smooth and strongly convex, $g$ CCP, $X_{k} \subseteq \mathbf{R}^n$;

\item Dual ascent (for the problem $\min_{x\in\mathbf{R}^n} f(x; t)$ subject to $A x = b$), for $f(\cdot; t_{k})$ strongly smooth and strongly convex;

\item D-R splitting, ADMM, doubly-regularized saddle-points, $\ldots$, for similar assumptions.

\end{itemize}

\subsection{Beyond strong convexity/strong smoothness}

We briefly touch here (and in this subsection alone) the more general case of relaxing the Assumption~\ref{as.1}, to generic convex problems. In particular, the previous tracking results can be extended also in case of more general $f(x;t)$, by using fixed-point theory in compact sets $X(t)$. 

E.g., for the projected gradient, if the function is only strongly smooth and $\alpha<2/L$, one can arrive at results of the form of 

\begin{itemize}
\item Average fixed-point residual tracking:
\begin{equation}
\frac{1}{T}\sum_{k=1}^T \|\Pi_{X_{k}}(x_{k} - \alpha \nabla_{x}f(x_k; t_{k})) - x_k\|^2 \leq O(1/T) + O(\tilde{K}D + \tilde{K}^2)
\end{equation}
where $\tilde{K}$ is a bound on a sequence of optimizers $\{x^*(t_{k})\}$, i.e., $\|x^*(t_{k+1}) - x^*(t_{k})\| \leq \tilde{K}$, (note that the optimizers need not be unique now); and $D = \max_t$ diam $X(t)$
\item Dynamic regret (aka objective function tracking):
\begin{equation}
\frac{1}{T}\underbrace{\sum_{k=1}^T f(x_{k}; t_{k}) - f(x_{k}^*; t_{k})}_{=:\, {\bf Reg}_T} \leq O(1/T) + O(\tilde{K}D + \tilde{K}^2) + O(K')
\end{equation}
where $\tilde{K}, D$ are as before, and $K'$ is a bound on the functional variations as $|f(x; t_{k+1}) - f(x; t_{k})| \leq {K}', \forall x \in X(t)$.
\end{itemize}

Similar results hold for other methods.

\section{Interlude: functions Lipschitz in time}

An interesting and useful result can be derived when the dependence of the cost function $f(x;t)$ over time is bounded in some sense. In particular, assume that $f$ has a well-defined gradient and the time derivative of the gradient in bounded, i.e.

\begin{assumption}\label{as.3}
The time derivative of the gradient of $f(x;t)$ in bounded uniformly in time,
$$
\|\nabla_{tx} f(x; t)\| \leq C_0, \quad \forall x\in\mathbf{R}^n, t\geq 0.
$$ 
\end{assumption}

Assumption~\ref{as.3} is more restrictive than Assumption~\ref{as.1} and it implies it
$$
\textrm{Assumption~\ref{as.3}} \implies \textrm{Assumption~\ref{as.1}},
$$ 
with $K = h\, C_0/m$, where we remind that $h$ is the sampling period. (One can see this in e.g.,~\cite{Dontchev2009}). 

Assumption~\ref{as.3} is a sort of Lipschitz condition in time, and when it is valid implies that all the running methods yield an asymptotical error of the order of $O(h)$ (i.e., linear in the sampling period).

\section{Prediction-correction algorithms}

Prediction-correction are methods that attempt at reducing the asymptotical error below $O(h)$. We look here at the results presented in~\cite{Paper1,Paper2,Paper3,Paper4}. 

To get better bounds, one needs stronger assumptions. In addition to Assumption~\ref{as.1}, here we will assume

\begin{assumption}\label{as.4}
Higher derivatives of the cost function are bounded as
$$
\|\nabla_{xxx} f(x; t)\| \leq C_1, \quad \|\nabla_{txx} f(x; t)\| \leq C_2, \quad \|\nabla_{ttx} f(x; t)\| \leq C_3
$$
uniformly in time and for all $x \in \mathbf{R}^n$.
\end{assumption}

Assumption~\ref{as.4} is an extension of Newton's assumptions (for Newton's method one requires $C_1$) that also requires the time variations of the Hessian and gradient to be bounded. 

Prediction-correction methods attempt at inferring how the optimizers are changing in time, by applying a pertinent Taylor's expansion of the optimality conditions. 

For example, if we were to solve the unconstrained problem
\begin{equation}
\min_{x \in \mathbf{R}^n} f(x; t),
\end{equation}
and we wanted to predict the optimizer at time $t_{k+1}$, only from data available at time $t_k$, one could start from the optimality condition at time $t_{k+1}$
\begin{equation}
\nabla_x f(x^*(t_{k+1}); t_{k+1}) = 0
\end{equation}
(which we can't solve) and Taylor expand as
\begin{equation}
\nabla_x f(x^*(t_{k+1}); t_{k+1}) \approx \nabla_x f(x^*(t_{k}); t_{k}) + h\, \nabla_{tx} f(x^*(t_{k}); t_{k}) + \nabla_{xx} f(x^*(t_{k}); t_{k}) \delta x = 0
\end{equation}
Since we don't have $x^*(t_k)$, we can substitute $x_k$ and obtain a class of prediction schemes:
\begin{equation}
\nabla_x f(x_k; t_{k}) + h\, \nabla_{tx} f(x_k; t_{k}) + \nabla_{xx} f(x_k; t_{k}) \delta x = \gamma \nabla_x f(x_k; t_{k}),
\end{equation}
and the prediction $x_{k+1|k}$ is given by
\begin{equation}
x_{k+1|k} = x_k + \delta x = x_k - [\nabla_{xx} f(x_k; t_{k})]^{-1}( h\, \nabla_{tx} f(x_k; t_{k}) +(1-\gamma) \nabla_{x} f(x_k; t_{k}) ).
\end{equation}
Here $\gamma \in [0,1]$ is an extra tuning parameter. For $\gamma = 1$, we have a tangential update: we are moving along the ``iso-suboptimal manifold''. When $\gamma = 0$, we have a Newton-like update, so that on top of predicting we are also going towards the optimizer. See a nice figure in~\cite{Simonetto2017}. 

If we were to solve the constrained problem
\begin{equation}
\min_{x \in X} f(x; t),
\end{equation}
then we would do prediction over the generalized inequality
\begin{equation}
\nabla_{x} f(x^*(t_{k+1}); t_{k+1}) + N_X(x^*(t_{k+1})) \ni 0,
\end{equation}
where $N_X$ is the normal cone operator, which leads to 
\begin{equation}
\nabla_{x} f(x_k; t_{k}) + h\, \nabla_{tx} f(x_k; t_{k}) + \nabla_{xx} f(x_k; t_{k})\,(x_{k+1|k}-x_k) +N_X(x_{k+1|k}) \ni {0}
\end{equation}
or equivalently, calling $Q_k = \nabla_{xx} f(x_k; t_{k})$, $c_k = h\, \nabla_{tx} f(x_k; t_{k})$,
\begin{equation}\label{qp}
x_{k+1|k} = \textrm{arg}\min_{y\in X} \{ 1/2 y^T Q_k y + c_k^T y\}. 
\end{equation}

Now, since we don't want to solve an optimization problem with another optimization problem (however easy), we can set up an approximate scheme for~\eqref{qp} as
\begin{equation}
y_{i+1} = \Pi_{X}[ y_i - \beta (Q_k y_i + c_k)]
\end{equation}
that is a projected gradient method that has to run for $P$ prediction steps and with stepsize $\beta>0$.

If we were to solve a linearly constrained problem, a similar construct would apply for both primal and dual variable in a dual ascent setting. 

\subsection{Prototypical algorithm}

As for the running methods, we report here a prototypical prediction-correction algorithm, here focussed on the projected gradient (but similar for gradient and dual ascent) 

\begin{itemize}
\item Time $t_0$, guess $x_0$
\item Time $t_{k}$
\begin{enumerate}
\item Set $Q_k = \nabla_{xx} f(x_k; t_{k})$, $c_k = h\, \nabla_{tx} f(x_k; t_{k})$
\item Set $y_{0} = x_{k}$
\item Perform $P$ prediction steps:
\begin{equation}\label{pred}
y_{i+1} = \Pi_{X}[ y_i - \beta (Q_k y_i + c_k)]
\end{equation}
\item Set $\tilde{x}_{k+1|k} = y_{P}$ (approximate prediction)
\end{enumerate}
\item Time $t_{k+1}$
\begin{enumerate}
\item Acquire a new function $f(\cdot; t_{k+1})$
\item Set $y_{0} = \tilde{x}_{k+1|k}$
\item Perform $C$ correction steps:
\begin{equation}\label{run}
y_{i+1} = \Pi_{X}[y_i - \alpha \nabla_{x} f(y_i; t_{k+1})]
\end{equation}
\item Set $x_{k+1} = y_{C}$
\end{enumerate}
\end{itemize}

Note that updates~\eqref{pred} are computationally cheap to carry out, once $Q_k$ and $c_k$ have been computed once, while updates~\eqref{run} may be more expensive. 

\subsection{Theoretical results}

Typical theoretical results goes as follows. Under Assumptions~\ref{as.1}, \ref{as.3}, \ref{as.4}, and a proper selection of stepsizes, number of prediction and correction steps, and sampling period, one is expect to track the solution trajectory up to a bound that depends on the problem properties and the sampling period. Depending on the method and on $P$ and $C$ we can have asymptotical errors that range from $O(h)$ to $O(h^4)$. 

We report the result for projected gradient in prediction-correction mode as defined in the previous subsection. 

\begin{theorem}\label{th.2}\emph{(Informal)}
The projected gradient method in prediction-correction mode generates a sequence $\{x_k\}$ as follows. 

Choose $\alpha,\beta < 2/L$. 

Under Assumptions~\ref{as.1}, \ref{as.3}, there exists a minimal number of prediction and correction steps $P, C$ for which globally
$$
\limsup_{k \to \infty} \|x_{k} - x^*(t_{k})\| = O(\varrho_1^C\, h)
$$
In addition, under Assumptions~\ref{as.4}, then locally (and for small $h$), there exists a minimal number of prediction and correction steps $P, C$ so that
$$
\limsup_{k \to \infty} \|x_{k} - x^*(t_{k})\| = \underbrace{O(\varrho_1^C \,h^2)}_{\textrm{prediction gain}} + \underbrace{O(\varrho_1^C \varrho_2^P\, h)}_{\textrm{approximation error}}
$$
where $\varrho_1,\varrho_2 < 1$, and $\varrho_1,\varrho_2$ are the contraction rates for $\alpha$ and $\beta$, respectively. 

Convergence is Q-linear in both cases.
\end{theorem} 

Theorem~\ref{th.2} says that tracking is not worse than correction-only method in the worst case. If the function has extra properties and we are interested in a local result, then a better asymptotical error can be achieved, provided some (stricter) conditions on the number of prediction and correction steps are verified. 

The asymptotical error is composed of two terms; one which is labeled as approximation error, which is due to the early termination of the prediction step (if $P \to \infty$ and prediction is exact, this term goes to $0$). The other, named prediction gain is the gain coming from using a prediction step, which brings the error down to a $O(h^2)$ dependence on the sampling period. 

If $C$ grows, then the error reduces, as expected. 

Theorem~\ref{th.2} can be modified for gradient methods and dual ascent methods. 

\section{A summary}

We give now a short comparison between running and prediction-correction methods. 

\begin{table}[H]
\caption{Comparison between correction-only and prediction-correction methods}
\centering
\begin{tabular}{ccc}
 & \textbf{Correction-only}  &    \textbf{Prediction-correction}  \\ \hline
Assumptions & Weak (mainly standard) & Stronger \\
Complexity & Low & Higher \\
Error & $O(h)$ & $O(h)$ - $O(h^4)$ \\
Methods & Many, see \cite{Simonetto20XX} & A few, see \cite{Paper1,Paper2,Paper3,Paper4} \\ \hline 
\end{tabular}
\end{table}

As one can see, correction-only method can tackle a larger class of problems up to a limited accuracy. Prediction-correction methods can achieve a better asymptotical error at the price of stronger assumptions and computational complexity. 

We note that, in some cases, even keeping the computational time fixed, prediction-correction may achieve better errors than correction-only methods. This is because prediction steps are computationally easier than correction steps, and one can trade-off a few correction steps for many prediction ones. So, prediction-correction are very relevant even in practice.

\section{Applications}

Many applications entail some degree of time-varying optimization. We report below a collection of tested applications (either in correction-only mode or prediction-correction). 

\begin{itemize}
\item Energy: e.g., time-varying optimal power flow and related, see~\cite{DallAnese2016,DallAnese2017,DallAnese2017a,Hauswirth2017,Tang2017,Liu2017,Zhang2017,Liu2018,Zhou2018}

\item Transportation: \cite{Su2009,Eser2018}

\item Robotics: e.g., dynamic consensus \cite{verscheure2009time,Zavlanos2013,Ling2013,Fazlyab2015,Fazlyab2016,Li2018}

\item Control: e.g., model predictive control~\cite{Jerez2014,Hours2014,Gutjahr2016,Anitescu2017, Paternain2018}

\item Signal processing: e.g., estimation in data streams, \cite{Asif2014,Yang2015,Vaswani2015,Balavoine2015,Simonetto2015a,Sopasakis2016,Maros2017}

\item Others: economics~\cite{Dontchev2013}, computational history~\cite{Nagurney2006}
\end{itemize}

\bibliographystyle{ieeetr}
\bibliography{PaperCollection00}

\begin{thebibliography}{10}

\bibitem{Polyak1987}
B.~T. Polyak, {\em Introduction to Optimization}.
\newblock Optimization Software, Inc., 1987.

\bibitem{Ye2015}
M.~Ye and G.~Hu, ``{Distributed Optimization for Systems with Time-Varying
  Quadratic Objective Functions},'' in {\em Proceedings of the 54th IEEE
  Conference on Decision and Control}, (Osaka, Japan), pp.~3285 -- 3290,
  December 2015.

\bibitem{Rahili2015a}
S.~Rahili, W.~Ren, and P.~Lin, ``{Distributed Convex Optimization of
  Time-Varying Cost Functions for Double-integrator Systems Using Nonsmooth
  Algorithms},'' in {\em Proceedings of the American Control Conference},
  (Chicago (IL), USA), pp.~68 -- 73, July 2015.

\bibitem{Rahili2015}
S.~Rahili and W.~Ren, ``{Distributed Convex Optimization for Continuous-Time
  Dynamics with Time-Varying Cost Functions},'' {\em IEEE Transactions on
  Automatic Control}, vol.~62, no.~4, pp.~1590 -- 1605, 2017.

\bibitem{Gong2016}
P.~Gong, F.~Chen, and W.~Lan, ``{Time-varying convex optimization for
  double-integrator dynamics over a directed network},'' in {\em Proceedings of
  the 35th Chinese Control Conference}, (Chengdu, China), pp.~7341 -- 7346,
  July 2016.

\bibitem{Fazlyab2015}
M.~Fazlyab, S.~Paternain, V.~Preciado, and A.~Ribeiro, ``{Interior Point Method
  for Dynamic Constrained Optimization in Continuous Time},'' in {\em
  Proceedings of the American Control Conference}, (Boston (MA), USA), pp.~5612
  -- 5618, July 2016.

\bibitem{Fazlyab2016}
M.~Fazlyab, S.~Paternain, V.~Preciado, and A.~Ribeiro, ``{Prediction-Correction
  Interior-Point Method for Time-Varying Convex Optimization},'' {\em IEEE
  Transactions on Automatic Control}, vol.~63, no.~7, 2018.

\bibitem{Moreau1977}
J.~J. Moreau, ``{Evolution Problem Associated with a Moving Convex Set in a
  Hilbert Space},'' {\em Journal of Differential Equations}, vol.~26, pp.~347
  -- 374, 1977.

\bibitem{Popkov2005}
A.~Y. Popkov, ``{Gradient Methods for Nonstationary Unconstrained Optimization
  Problems},'' {\em Automation and Remote Control}, vol.~66, no.~6, pp.~883 --
  891, 2005.
\newblock Translated from Avtomatika i Telemekhanika, No. 6, 2005, pp. 38 --
  46.

\bibitem{Tu2011}
S.-Y. Tu and A.~H. Sayed, ``{Mobile Adaptive Networks},'' {\em IEEE Journal of
  Selected Topics in Signal Processing}, vol.~5, no.~4, pp.~649 -- 664, 2011.

\bibitem{Bajovic2011}
D.~Bajovic, D.~Jakovetic, J.~Xavier, B.~Sinopoli, and J.~M.~F. Moura,
  ``{Distributed Detection via Gaussian Running Consensus: Large Deviations
  Asymptotic Analysis},'' {\em IEEE Transactions on Signal Processing},
  vol.~59, no.~9, pp.~4381 -- 4396, 2011.

\bibitem{Dontchev2013}
A.~L. Dontchev, M.~I. Krastanov, R.~T. Rockafellar, and V.~M. Veliov, ``{An
  Euler-Newton Continuation method for Tracking Solution Trajectories of
  Parametric Variational Inequalities},'' {\em SIAM Journal of Control and
  Optimization}, vol.~51, no.~51, pp.~1823 -- 1840, 2013.

\bibitem{Zavlanos2013}
M.~M. Zavlanos, A.~Ribeiro, and G.~J. Pappas, ``{Network Integrity in Mobile
  Robotic Networks},'' {\em IEEE Transactions on Automatic Control}, vol.~58,
  no.~1, pp.~3 -- 18, 2013.

\bibitem{Jakubiec2013}
F.~Y. Jakubiec and A.~Ribeiro, ``{D-MAP: Distributed Maximum a Posteriori
  Probability Estimation of Dynamic Systems},'' {\em IEEE Transactions on
  Signal Processing}, vol.~61, no.~2, pp.~450 -- 466, 2013.

\bibitem{Ling2013}
Q.~Ling and A.~Ribeiro, ``{Decentralized Dynamic Optimization Through the
  Alternating Direction Method of Multipliers},'' {\em IEEE Transactions on
  Signal Processing}, vol.~62, no.~5, pp.~1185 -- 1197, 2014.

\bibitem{Simonetto2014c}
A.~Simonetto and G.~Leus, ``{Distributed Asynchronous Time-Varying Constrained
  Optimization},'' in {\em Proceedings of the Asilomar Conference on Signals,
  Systems, and Computers}, (Pacific Grove, USA), November 2014.

\bibitem{Simonetto2014d}
A.~Simonetto and G.~Leus, ``{Double Smoothing for Time-Varying Distributed
  Multi-user Optimization},'' in {\em Proceedings of the IEEE Global Conference
  on Signal and Information Processing}, (Atlanta, US), December 2014.

\bibitem{Xi2016a}
C.~Xi and U.~A. Khan, ``{Distributed Dynamic Optimization over Directed
  Graphs},'' in {\em Proceedings of the 55th IEEE Conference on Decision and
  Control}, (Las Vegas, NV, US), pp.~245 -- 250, December 2016.

\bibitem{Sun2017}
C.~Sun, M.~Ye, and G.~Hu, ``{Distributed Time-varying Quadratic Optimization
  for Multiple Agents under Undirected Graphs},'' {\em IEEE Transactions on
  Automatic Control}, vol.~62, no.~7, pp.~3687 -- 3694, 2017.

\bibitem{Maros2017}
M.~Maros and J.~Jalden, ``{ADMM for Distributed Dynamic Beam-forming},'' {\em
  IEEE Transactions on Signal and Information Processing over Networks},
  vol.~4, no.~2, pp.~220 -- 235, 2018.

\bibitem{Simonetto20XX}
A.~Simonetto, ``{Time-Varying Convex Optimization via Time-Varying Averaged
  Operators },'' {\em arXiv: 1704.07338v1}, 2017.

\bibitem{Robinson1980}
S.~M. Robinson, ``{Strongly Regular Generalized Equations},'' {\em Mathematics
  of Operations Research}, vol.~5, no.~1, pp.~43 -- 62, 1980.

\bibitem{Dontchev2009}
A.~L. Dontchev and R.~T. Rockafellar, {\em Implicit Functions and Solution
  Mappings}.
\newblock Springer, 2009.

\bibitem{Zavala2010}
V.~M. Zavala and M.~Anitescu, ``{Real-Time Nonlinear Optimization as a
  Generalized Equation},'' {\em SIAM Journal of Control and Optimization},
  vol.~48, no.~8, pp.~5444 -- 5467, 2010.

\bibitem{Kungurtsev2017}
V.~Kungurtsev and J.~J{\"a}schke, ``{A Prediction-Correction Path-Following
  Algorithm for Dual-Degenerate Parametric Optimization Problems},'' {\em SIAM
  Journal on Optimization}, vol.~27, no.~1, pp.~538 -- 564, 2017.

\bibitem{Allgower1990}
E.~L. Allgower and K.~Georg, {\em Numerical Continuation Methods: An
  Introduction}.
\newblock Springer-Verlag, 1990.

\bibitem{Cojocaru2005}
M.~G. Cojocaru, P.~Daniele, and A.~Nagurney, ``{Projected Dynamical Systems and
  Evolutionary Variational Inequalities via Hilbert Spaces with
  Applications},'' {\em Journal of Optimization Theory and Applications},
  vol.~127, no.~3, pp.~549 -- 563, 2005.

\bibitem{Nagurney2006}
A.~Nagurney and J.~Pan, ``{Evolution Variational Inequalities and Projected
  Dynamical Systems with Application to Human Migration},'' {\em Mathematical
  and Computer Modelling}, vol.~43, no.~5 -- 6, pp.~646 -- 657, 2006.

\bibitem{Nesterov2012}
Y.~Nesterov, ``{Towards Non-symmetric Conic Optimization},'' {\em Optimization
  Methods and Software}, vol.~27, no.~4 -- 5, pp.~893 -- 917, 2012.

\bibitem{Paper1}
A.~Simonetto, A.~Mokhtari, A.~Koppel, G.~Leus, and A.~Ribeiro, ``{A Class of
  Prediction-Correction Methods for Time-Varying Convex Optimization},'' {\em
  IEEE Transactions on Signal Processing}, vol.~64, no.~17, pp.~4576 -- 4591,
  2016.

\bibitem{Paper2}
A.~Simonetto, A.~Koppel, A.~Mokhtari, G.~Leus, and A.~Ribeiro, ``{Decentralized
  Prediction-Correction Methods for Networked Time-Varying Convex
  Optimization},'' {\em IEEE Transactions on Automatic Control}, vol.~62,
  no.~11, pp.~5724 -- 5738, 2017.

\bibitem{Paper3}
A.~Simonetto and E.~{Dall'Anese}, ``{Prediction-Correction Algorithms for
  Time-Varying Constrained Optimization},'' {\em IEEE Transactions on Signal
  Processing}, vol.~65, no.~20, pp.~5481 -- 5494, 2017.

\bibitem{Paper4}
A.~Simonetto, ``{Dual Prediction-Correction Methods for Linearly Constrained
  Time-Varying Convex Programs},'' {\em arXiv: 1709.05850}, 2017.

\bibitem{Guddat1990}
J.~Guddat and F.~{Guerra Vazquez and H.~T. Jongen}, {\em {Parametric
  Optimization: Singularities, Pathfollowing and Jumps}}.
\newblock John Wiley \& Sons, Chichester, UK, 1990.

\bibitem{Simonetto2017}
A.~Simonetto and E.~{Dall'Anese}, ``{A First -Order Prediction-Correction
  Algorithm for Time-Varying (Constrained) Optimization},'' in {\em Proceedings
  of the 20th IFAC World Congress}, (Toulouse, France), July 2017.

\bibitem{DallAnese2016}
E.~Dall'Anese and A.~Simonetto, ``{Optimal Power Flow Pursuit},'' {\em IEEE
  Transactions on Smart Grid}, vol.~9, no.~2, pp.~942 -- 952, 2018.

\bibitem{DallAnese2017}
E.~{Dall'Anese}, S.~Guggilam, A.~Simonetto, Y.~C. Chen, and S.~V. Dhople,
  ``{Optimal Regulation of Virtual Power Plants},'' {\em IEEE Transactions on
  Power Systems}, vol.~33, no.~2, pp.~1868 -- 1881, 2018.

\bibitem{DallAnese2017a}
E.~{Dall'Anese}, A.~Bernstein, and A.~Simonetto, ``{Feedback-based
  projected-gradient method for real-time optimization of aggregations of
  energy resources},'' in {\em Proceedings of the Global Conference on Signal
  and Information Processing}, (Montreal, QC, Canada), November 2017.

\bibitem{Hauswirth2017}
A.~Hauswirth, A.~Zanardi, S.~Bolognani, F.~D\"orfler, and G.~Hug, ``{Online
  optimization in closed loop on the power flow manifold},'' in {\em
  Proceedings of the IEEE PowerTech conference}, (Manchester, UK), June 2017.

\bibitem{Tang2017}
Y.~Tang, K.~Dvijotham, and S.~Low, ``{Real-Time Optimal Power Flow},'' {\em
  IEEE Transactions on Smart Grid}, vol.~8, no.~6, pp.~2963 -- 2973, 2017.

\bibitem{Liu2017}
H.~J. Liu, W.~Shi, and H.~Zhu, ``{Decentralized Dynamic Optimization for Power
  Network Voltage Control},'' {\em IEEE Transactions on Signal and Information
  Processing over Networks}, vol.~3, no.~3, pp.~568 -- 579, 2017.

\bibitem{Zhang2017}
Y.~Zhang, E.~{Dall'Anese}, and M.~Hong, ``{Dynamic ADMM for real-time optimal
  power flow},'' in {\em Proceedings of the IEEE Global Conference on Signal
  and Information Processing}, (Montreal, QC, Canada), November 2017.

\bibitem{Liu2018}
J.~Liu, J.~Marecek, A.~Simonetto, and M.~Takac, ``{A Coordinate-Descent
  Algorithm for Tracking Solutions in Time-Varying Optimal Power Flows},'' in
  {\em Proceedings of the XX Power Systems Computation Conference}, (Dublin,
  Ireland), June 2018.

\bibitem{Zhou2018}
X.~Zhou, E.~{Dall'Anese}, L.~Chen, and A.~Simonetto, ``{An Incentive-Based
  Online Optimization Framework for Distribution Grids},'' {\em IEEE
  Transactions on Automatic Control}, vol.~63, no.~7, 2018.

\bibitem{Su2009}
W.~Su, {\em {Traffic Engineering and Time-Varying Convex Optimization}}.
\newblock PhD thesis, The Pennsylvania State University, May 2009.

\bibitem{Eser2018}
E.~Eser, J.~Monteil, and A.~Simonetto, ``{On the Tracking of Dynamical Optimal
  Meeting Points},'' in {\em Proceedings of the 15th IFAC Symposium on Control
  in Trasportation Systems}, (Savona, Italy), June 2018.

\bibitem{verscheure2009time}
D.~Verscheure, B.~Demeulenaere, J.~Swevers, J.~De~Schutter, and M.~Diehl,
  ``{Time-Optimal Path Tracking for Robots: a Convex Optimization Approach},''
  {\em IEEE Transactions on Automatic Control}, vol.~54, no.~10, pp.~2318 --
  2327, 2009.

\bibitem{Li2018}
J.~Li, M.~Mao, F.~Uhlig, and Y.~Zhang, ``{Z-type neural-dynamics for
  time-varying nonlinear optimization under a linear equality constraint with
  robot application},'' {\em Journal of Computational and Applied Mathematics},
  vol.~327, no.~1, pp.~155 -- 166, 2018.

\bibitem{Jerez2014}
J.~L. Jerez, P.~J. Goulart, S.~Richter, G.~A. Constantinides, E.~C. Kerrigan,
  and M.~Morari, ``{Embedded Online Optimization for Model Predictive Control
  at Megahertz Rates},'' {\em IEEE Transactions on Automatic Control}, vol.~59,
  no.~12, pp.~3238 -- 3251, 2014.

\bibitem{Hours2014}
J.-H. Hours and C.~N. Jones, ``{A Parametric Non-Convex Decomposition Algorithm
  for Real-Time and Distributed NMPC},'' {\em IEEE Transactions on Automatic
  Control}, vol.~61, no.~2, pp.~287 -- 302, 2016.

\bibitem{Gutjahr2016}
B.~Gutjahr, L.~Gr{\"o}ll, and M.~Werling, ``{Lateral Vehicle Trajectory
  Optimization Using Constrained Linear Time-Varying MPC},'' {\em IEEE
  Transactions on Intelligent Transportation Systems}, pp.~1 -- 10, 2016.

\bibitem{Anitescu2017}
M.~Anitescu and V.~Zavala, ``{MPC as a DVI: Implications on Sampling Rates and
  Accuracy},'' in {\em Proceedings of the Control and Decision Conference},
  (Melbourne, Australia), December 2017.

\bibitem{Paternain2018}
S.~Paternain, M.~Morari, and A.~Ribeiro, ``{A Prediction-Correction Method for
  Model Predictive Control},'' in {\em Proceedings of the American Control
  Conference}, (Milwaukee, WI, USA), June 2018.

\bibitem{Asif2014}
M.~S. Asif and J.~Romberg, ``{Sparse recovery of streaming signals using
  $\ell_1$-homotopy },'' {\em IEEE Transactions on Signal Processing}, vol.~62,
  no.~16, pp.~4209 -- 4223, 2014.

\bibitem{Yang2015}
Y.~Yang, M.~Zhang, M.~Pesavento, and D.~P. Palomar, ``{An Online Parallel and
  Distributed Algorithm for Recursive Estimation of Sparse Signals},'' {\em
  IEEE Transactions on Signal and Information Processing over Networks},
  vol.~2, no.~3, pp.~290 -- 305, 2016.

\bibitem{Vaswani2015}
N.~Vaswani and J.~Zhan, ``{Recursive Recovery of Sparse Signal Sequences from
  Compressive Measurements: A Review},'' {\em IEEE Transactions on Signal
  Processing}, vol.~64, no.~13, pp.~3523 -- 3549, 2016.

\bibitem{Balavoine2015}
A.~Balavoine, J.~Romberg, and C.~Rozell, ``{Discrete and continuous iterative
  soft thresholding with a dynamic input},'' {\em IEEE Transactions on Signal
  Processing}, vol.~63, no.~12, pp.~3165 -- 3176, 2015.

\bibitem{Simonetto2015a}
A.~Simonetto and G.~Leus, ``{On Non-Differentiable Time-Varying
  Optimization},'' in {\em Proceedings of the 6th IEEE CAMSAP}, (Cancun,
  Mexico), December 2015.

\bibitem{Sopasakis2016}
P.~Sopasakis, N.~Freris, and P.~Patrinos, ``{Accelerated Reconstruction of a
  Compressively Sampled Data Stream},'' in {\em Proceedings of the 24th
  EUSIPCO}, (Budapest, Hungary), pp.~1078 -- 1082, September 2016.

\end{thebibliography}
\end{document}